\documentclass[12pt]{amsart}
\usepackage{amsmath, amssymb, amsthm}
\usepackage{hyperref}
\usepackage{geometry}
\newtheorem{claim}{Claim}
\newtheorem{question}{Question}

\title{On the Maximum Volume Solid Wrappable by a Given Sheet of Paper}
\author{R.~Nandakumar}
\email{nandacumar@gmail.com}
\date{\today}

\begin{document}

\maketitle

\begin{abstract}
We consider the problem of wrapping three-dimensional solid bodies with a
given planar sheet of paper, where the paper may be folded or wrinkled
but not stretched or torn. We propose a conjecture characterising the
maximum-volume solid wrappable by any given sheet: the maximum is always
achieved (or approached) by a \emph{non-convex} body. In other words,
for any convex solid wrappable by a given sheet, there exists a
non-convex solid of strictly greater volume that the same sheet can wrap.
We discuss related work, a key subquestion involving the sphere, and
several further directions.
\end{abstract}

\section{Setup and definitions}

Let $D \subset \mathbb{R}^2$ be a connected planar region --- a
\emph{sheet of paper} --- that is non-degenerate in the sense that it
has positive area and no infinitely thin parts. The sheet may be
polygonal or more general.

We say that $D$ \emph{wraps} a compact connected solid $B \subset
\mathbb{R}^3$ if the sheet can be placed so as to enclose $B$ in the
following topological sense: every path from any point on the surface
$\partial B$ to any point at infinity in $\mathbb{R}^3$ must cross the
image of $D$ at least once. The paper may be folded or wrinkled
arbitrarily, but may not be stretched or torn.

Formally, following Karasev \cite{karasev}, this is the question of
bounding the maximal volume enclosed by a $1$-Lipschitz image of $D$
in $\mathbb{R}^3$, where the boundary $\partial D$ is required to be
glued to itself in the image (so that the image forms a closed surface).

We call a wrap \emph{neat} if every point on $\partial B$ is in contact
with exactly one point of $D$, so that the surface area of $B$ equals
the area of $D$ exactly. For a wrap maximising the volume of $B$, one
expects the surface area of $B$ to be \emph{strictly less} than the
area of $D$, with the excess paper folded or crumpled.

\section{The main conjecture}

\begin{claim}
Let $D \subset \mathbb{R}^2$ be any connected, non-degenerate planar
region. Among all connected compact solids $B \subset \mathbb{R}^3$
that can be wrapped by $D$, the supremum of $\mathrm{Vol}(B)$ is never
achieved by a convex body. Equivalently, for every convex solid $B$
wrappable by $D$, there exists a non-convex solid $B'$, also wrappable
by $D$, with $\mathrm{Vol}(B') > \mathrm{Vol}(B)$.
\end{claim}

A constrained version of the conjecture, in which $D$ is itself
required to be convex, may be easier to approach and is also of
interest.

\section{Relation to existing work}

Work of Bleecker \cite{bleecker} and Pak \cite{pak} on isometric
foldings and inflation of polyhedral surfaces appears closely related.
Their results pertain primarily to neat wraps, where the paper's metric
is preserved. It is conceivable that the present conjecture is a
corollary of their work, but the non-neat case --- where the paper is
allowed to crumple --- appears to go beyond current machinery.

Karasev \cite{karasev} observed that the reverse question --- given a
convex body $B$, what is the minimum disk that wraps it? --- has a
clean answer: a disk of radius equal to the intrinsic diameter of
$\partial B$, attained by the exponential map. This uses the fact that
$\partial B$ for a convex body has Alexandrov comparison curvature
$\geq 0$, which forces exponential maps to be $1$-Lipschitz. This
resolves the reverse question but does not directly address the forward
problem of maximising volume for a given sheet.

\section{A key subquestion: wrapping a sphere with a disk}

The exponential map from a point on the unit sphere $S^2$ produces a
disk of radius $\pi$ (area $\pi^3 \approx 31.0$) that maps surjectively
onto $S^2$ (surface area $4\pi \approx 12.6$), a ratio of approximately
$2.47$. This disk wraps the unit sphere non-neatly, with substantial
excess paper.

\begin{question}
Can a disk of radius $\pi$ wrap a non-convex solid of volume strictly
greater than $\frac{4}{3}\pi$ (the volume of the unit sphere)?
\end{question}

The intuition is that the large excess of paper --- a factor of $2.47$
in area --- should allow coverage of a solid with greater volume than
the sphere, provided the excess can be folded into concavities of a
non-convex solid in a topologically valid way. Making this precise
appears non-trivial: the exact manner in which the excess paper must be
folded and creased requires careful analysis, and one may need a disk of
radius strictly greater than $\pi$ to handle the boundary behaviour.

A positive answer to this subquestion would constitute the first
concrete instance of the main conjecture, for the specific case of a
disk sheet.

\section{Further questions}

\begin{enumerate}
\item \textbf{Singularities.} If an optimal (or supremum-approaching)
solid exists, what is the nature of the singularities on its surface?
The non-neat setting allows for highly irregular boundaries, and
characterising the possible singularities would shed light on the
geometry of the problem.

\item \textbf{Multi-layered wrapping.} One can define a \emph{$k$-layered
wrap} by requiring that every path from infinity to any point on
$\partial B$ crosses the image of $D$ at least $k$ times. What is the
maximum volume solid that is $2$-layered wrappable by a given sheet?
The additional constraint may make the problem more tractable.

\item \textbf{Attainment.} Is the supremum of wrappable volumes always
attained, or only approached by a sequence of solids? If the supremum
is never attained, this itself implies that convex bodies are never
optimal, provided one can always find non-convex bodies approaching the
supremum.
\end{enumerate}

\section*{Acknowledgements}

The author thanks Roman Karasev for comments and reformulations
communicated privately \cite{karasev}, and Joseph O'Rourke for comments
on the MathOverflow post \cite{orourke}, both of which clarified the
structure of the problem considerably.


\begin{thebibliography}{9}

\bibitem{bleecker}
D.~Bleecker,
\emph{Volume increasing isometric deformations of convex polyhedra},
J.\ Differential Geom.\ \textbf{43} (1996), no.~3, 505--526.

\bibitem{pak}
I.~Pak,
\emph{Inflating polyhedral surfaces},
Preprint, Department of Mathematics, MIT, no.~326, 2006.

\bibitem{mo}
R.~Nandakumar,
\emph{A claim on wrapping 3D solids with a sheet of paper},
MathOverflow, 2026.
\url{https://mathoverflow.net/questions/508997/a-claim-on-wrapping-3d-solids-with-a-sheet-of-paper}

\bibitem{karasev}
R.~Karasev,
personal communication, 2026.

\bibitem{orourke}
J.~O'Rourke,
comments on \cite{mo}, 2026.

\end{thebibliography}
\end{document}